\documentclass[a4paper,12pt,reqno]{amsart}
\usepackage{showkeys}
\usepackage{mathtools}
\usepackage{amsfonts,amsmath}
\usepackage{amssymb,amsthm}		
\usepackage{amsthm}
\usepackage[T1]{fontenc}
\usepackage{lmodern}
\usepackage{hyperref}
\usepackage[all]{hypcap}
\usepackage{fullpage}
\usepackage{color}
\usepackage{graphicx}
\usepackage[dvipsnames]{xcolor}
\usepackage{fancyhdr}
\usepackage{psfrag}

\usepackage{indentfirst}
\usepackage{float}
\usepackage{latexsym}
\usepackage{graphics}
\usepackage{verbatim}

\newtheorem{theorem}{Theorem}
\newtheorem*{theorem*}{Theorem}

\newtheorem{corollary}{Corollary}

\newtheorem{proposition}{Proposition}
\theoremstyle{definition}

\theoremstyle{remark}

\begin{document}

\author{Stephan Baier, Arpit Bansal}
\address{Stephan Baier, Ramakrishna Mission Vivekananda Educational Research Institute, Department of Mathematics, 
G. T. Road, PO Belur Math, Howrah, West Bengal 711202, 
India; email: email$_{-}$baier@yahoo.de}
\address{Arpit Bansal, Jawaharlal Nehru University, School of Physical Sciences, New-Delhi 110067, India; email: apabansal@gmail.com}
\title{The large sieve with power moduli for $\mathbb{Z}[i]$}
\begin{abstract} We establish a large sieve inequality for power moduli in $\mathbb{Z}[i]$, extending earlier work by L. Zhao and the first-named author on the large sieve for power moduli for the classical case of moduli in $\mathbb{Z}$.  Our method starts with a version of the
large sieve for $\mathbb{R}^2$. We convert the resulting counting problem back into one for $\mathbb{Z}[i]$ which we then attack using Weyl differencing and Poisson summation. 
\end{abstract}
\maketitle

\section{Introduction}
The classical large sieve inequality with additive characters asserts that
\begin{equation*} 
\sum\limits_{q\le Q} \sum\limits_{\substack{a=1\\ (a,q)=1}}^q \left| \sum\limits_{M<n\le M+N} a_n e\left(n\cdot \frac{a}{q}\right) \right|^2
\le \left(Q^2+N-1\right)\sum\limits_{M<n\le M+N} |a_n|^2,
\end{equation*}
where $Q,N\in \mathbb{N}$ and $M\in \mathbb{Z}$. This inequality has numerous applications in analytic number theory, in  particular, in sieve 
theory and to questions regarding the distribution of arithmetic functions in arithmetic progressions. 

The large sieve with resticted sets of 
moduli $q$, in particular power moduli, was considered in a series of papers by Baier, Zhao and Halupczok (see \cite{Bai}, \cite{BZ1},
\cite{BZ2}, \cite{Hal} and \cite{Zh1}), and these results turned out to be useful tools for 
applications (see \cite{BFKS} and \cite{BPS}, for example). In the case of square moduli, it was first established by Zhao \cite{Zh1} that 
\begin{equation} \label{firstls}
\begin{split}
& \sum\limits_{q\le Q} \sum\limits_{\substack{a=1\\ (a,q)=1}}^{q^2} \left| \sum\limits_{M<n\le M+N} a_n e\left(n\cdot \frac{a}{q^2}\right) 
\right|^2 \\
\ll_{\varepsilon} & (QN)^{\varepsilon}\left(Q^3+N\sqrt{Q}+\sqrt{N}Q^2\right)\sum\limits_{M<n\le M+N} |a_n|^2.
\end{split}
\end{equation}
This was improved in \cite{Bai}, where the term $N\sqrt{Q}$ on the right-hand side of \eqref{firstls} was replaced by $N$. A further 
improvement was obtained in \cite{BZ1}, where \eqref{firstls} with $N+\min\left\{N\sqrt{Q},\sqrt{N}Q^2\right\}$ in place of 
$N\sqrt{Q}+\sqrt{N}Q^2$ was established. In \cite{Zh1}, Zhao conjectured that the bound
\begin{equation} \label{conjec}
\sum\limits_{q\le Q} \sum\limits_{\substack{a=1\\ (a,q)=1}}^{q^k} \left| \sum\limits_{M<n\le M+N} a_n e\left(n\cdot \frac{a}{q^k}\right) \right|^2
\ll_{\varepsilon} (QN)^{\varepsilon}\left(N+Q^{k+1}\right)\sum\limits_{M<n\le M+N} |a_n|^2
\end{equation}
should hold for $k$-th power moduli ($k\in \mathbb{N}$ arbitrary but fixed). This conjecture is still open for every $k\ge 2$. 
In the same paper \cite{Zh1}, he established that
\begin{equation} \label{kls}
\begin{split}
& \sum\limits_{q\le Q} \sum\limits_{\substack{a=1\\ (a,q)=1}}^{q^k} \left| \sum\limits_{M<n\le M+N} a_n e\left(n\cdot \frac{a}{q^k}\right) 
\right|^2 \\
\ll_{\varepsilon} & (QN)^{\varepsilon}\left(Q^{k+1}+NQ^{1-1/\kappa}+N^{1-1/\kappa}Q^{1+k/\kappa}
\right) \sum\limits_{M<n\le M+N} |a_n|^2,
\end{split}
\end{equation}
where $\kappa=2^{k-1}$, thus generalizing \eqref{firstls}. Improvements of this result have been established in \cite{BZ2} and \cite{Hal}.

The large sieve for additive characters was extended to number fields by Huxley. In the case of the number field  $\mathbb{Q}(i)$ it takes the 
form
\begin{equation} \label{Huxley}
\sum\limits_{\substack{q\in \mathbb{Z}[i]\setminus\{0\}\\ \mathcal{N}(q)\le Q}} 
\sum\limits_{\substack{r \bmod{q}\\ (r,q)=1}} \left|\sum\limits_{\substack{n\in \mathbb{Z}[i]\\ \mathcal{N}(n)\le N}}  
a_n \cdot e\left(\Re\left(\frac{nr}{q}\right)\right)\right|^2 \ll \left(Q^2+N\right)\sum\limits_{\substack{n\in \mathbb{Z}[i]\\ 
\mathcal{N}(n)\le N}} |a_n|^2.
\end{equation}
Here as in the following $\mathcal{N}(q)$ denotes the norm of $q\in \mathbb{Z}[i]$, given by
$$
\mathcal{N}(q):=\Re(q)^2+\Im(q)^2.
$$
The large sieve with square norm moduli for the number field $\mathbb{Q}(i)$ was investigated in \cite{Ba2}, where an analogue 
of \eqref{firstls}
was established, namely the inequality
\begin{equation} \label{squarenorm}
\begin{split}
& \sum\limits_{\substack{q\in \mathbb{Z}[i]\setminus\{0\}\\ \mathcal{N}(q)\le Q^2\\ \mathcal{N}(q)=\Box}} \sum\limits_{\substack{r \bmod{q}\\ (r,q)=1}} \left|\sum\limits_{\substack{n\in \mathbb{Z}[i]\\ \mathcal{N}(n)\le N}}  
a_n \cdot e\left(\Re\left(\frac{nr}{q}\right)\right)\right|^2 \\
\ll & (QN)^{\varepsilon}\left(Q^3+Q^2\sqrt{N}+\sqrt{Q}N\right)\sum\limits_{\substack{n\in \mathbb{Z}[i]\\ \mathcal{N}(n)\le N}} |a_n|^2.
\end{split}
\end{equation}

In this paper, we go a step further and prove an analogue of \eqref{kls} for $\mathbb{Q}(i)$, i.e. a large sieve inequality with $k$-th power 
moduli for $\mathbb{Q}(i)$. 
Our approach will be more elegant than the previous one in \cite{Ba2}, where the double large sieve and lattice point counting in $\mathbb{R}^2$ 
were used. Here our method starts with a version of the
large sieve for $\mathbb{R}^2$. Then we convert the resulting counting problem back into one for $\mathbb{Z}[i]$ which can be attacked 
along similar lines as in \cite{Zh1} using Weyl differencing and Poisson summation. We begin with square moduli,
for which we obtain the essentially same bound as for square norm moduli in \eqref{squarenorm}. Then we generalize our method to $k$-th power 
moduli.

\section{Statement of main results}
We shall establish the following large sieve inequality for square moduli in $\mathbb{Z}[i]$. 

\begin{theorem} \label{squaremodZi} Let $Q,N\ge 1$ and $(a_n)_{n\in \mathbb{Z}[i]}$ be any sequence of complex numbers. Then 
\begin{equation*}
\begin{split}
& \sum\limits_{\substack{q\in \mathbb{Z}[i]\setminus\{0\}\\ \mathcal{N}(q)\le Q}} 
\sum\limits_{\substack{r \bmod{q^2}\\ (r,q)=1}} \left|\sum\limits_{\substack{n\in \mathbb{Z}[i]\\ \mathcal{N}(n)\le N}}  
a_n \cdot e\left(\Re\left(\frac{nr}{q^2}\right)\right)\right|^2\\ \ll & (QN)^{\varepsilon}\left(Q^3+Q^2\sqrt{N}+\sqrt{Q}N\right)
\sum\limits_{\substack{n\in \mathbb{Z}[i]\\ \mathcal{N}(n)\le N}} |a_n|^2,
\end{split}
\end{equation*}
where $\varepsilon$ is any positive constant, and the implied $\ll$-constant depends only on $\varepsilon$. 
\end{theorem}

Theorem \ref{squaremodZi} will then be generalized to $k$-th power moduli, for which we establish the following.

\begin{theorem} \label{powermodZi} Let $k\in \mathbb{N}$, $Q,N\ge 1$ and $(a_n)_{n\in \mathbb{Z}[i]}$ be any sequence of complex numbers.
Set $\kappa:=2^{k-1}$. Then
\begin{equation*}
\begin{split}
& \sum\limits_{\substack{q\in \mathbb{Z}[i]\setminus\{0\}\\ \mathcal{N}(q)\le Q}} 
\sum\limits_{\substack{r \bmod{q^k}\\ (r,q)=1}} \left|\sum\limits_{\substack{n\in \mathbb{Z}[i]\\ \mathcal{N}(n)\le N}}  
a_n \cdot e\left(\Re\left(\frac{nr}{q^k}\right)\right)\right|^2\\ \ll & 
(QN)^{\varepsilon}\left(Q^{k+1}+NQ^{1-1/\kappa}+N^{1-1/\kappa}Q^{1+k/\kappa}
\right)
\sum\limits_{\substack{n\in \mathbb{Z}[i]\\ \mathcal{N}(n)\le N}} |a_n|^2,
\end{split}
\end{equation*}
where $\varepsilon$ is any positive constant, and the implied $\ll$-constant depends only on $k$ and $\varepsilon$. 
\end{theorem}

\section{Large sieve for $\mathbb{R}^d$}
We shall employ the following version of the large sieve for $\mathbb{R}^d$ (in fact, we shall need it for the case $d=2$ only).

\begin{theorem} \label{ls} Let $R,N\in \mathbb{N}$, $N\ge 2$, $x_1,...,x_R\in \mathbb{R}^d$ and $(a_n)_{n\in \mathbb{Z}^d}$ 
be any $d$-fold sequence of complex numbers. Then 
$$
\sum\limits_{i=1}^R \left|\sum\limits_{\substack{n\in \mathbb{Z}^d\\ ||n||_2\le N^{1/d}}} 
a_{n} \cdot e\left(n\cdot x_i\right)\right|^2 \ll KNZ,
$$
where 
$$
K:=\max\limits_{1\le i\le R} \sharp \left\{j\in \{1,...,R\} : \min\limits_{z\in \mathbb{Z}^d} ||x_j-x_i-z||_2\le \sqrt{d}N^{-1/d}\right\} 
$$
and
\begin{equation} \label{Zdef}
Z:=\sum\limits_{\substack{n\in \mathbb{Z}^d\\ ||n||_2\le N^{1/d}}} |a_n|^2.
\end{equation}
\end{theorem}

Here as in the following, $||s||_2$ denotes the Euclidean norm of $s\in \mathbb{R}^d$, given by 
$$
||(s_1,s_2,..,s_d)||_2=\sqrt{s_1^2+s_2^2+\cdots + s_d^2}.  
$$

To prove Theorem \ref{ls}, we use the duality principle and the Poisson summation formula for $\mathbb{R}^d$.

\begin{proposition}[Duality principle, Theorem 288 in \cite{HLP}] \label{duality} Let $C = [c_{mn}]$ be a finite matrix with complex entries. 
The following two statements are equivalent: 
\begin{enumerate}
\item For any complex numbers $a_n$, we have 
\begin{align*}
\sum_m \mathrel \Big |\sum_n a_n c_{mn}\Big |^2 \leq \Delta \sum_n |a_n|^2.
\end{align*}
\item For any complex numbers $b_m$, we have
\begin{align*}
\sum_n \mathrel \Big |\sum_m b_m c_{mn}\Big |^2 \leq \Delta \sum_m |b_m|^2. 
\end{align*}
\end{enumerate}
\end{proposition}

\begin{proposition}[Poisson summation formula, see \cite{StW}] \label{poisson} 
Let $f:\mathbb{R}^d \rightarrow \mathbb{C}$ be a smooth function of rapid decay and $\Lambda$ be a lattice 
of full rank in $\mathbb{R}^d$. Then
$$
\sum\limits_{y\in \Lambda} f(y) = \frac{1}{\mbox{\rm Vol}(\mathbb{R}^d/\Lambda)} \cdot \sum\limits_{x\in \Lambda'} \hat{f}(x),
$$
where $\Lambda'$ is the dual lattice and $\hat{f}$ is the Fourier transform of $f$, defined as 
$$
\hat{f}(x)=\int\limits_{\mathbb{R}^2} f(y)e\left(-x\cdot y\right)\ dy.
$$
\end{proposition}

Here as in the following, by {\it rapid decay} we mean that the function $f:\mathbb{R}^d \rightarrow \mathbb{C}$ in question satisfies the bound
$$
f(y)\ll \left(1+||y||_2\right)^{-C}
$$
for some $C>d$. 

By a linear change of variables, we immediately deduce the following more general version of the Poisson summation formula for shifted 
lattices from Proposition \ref{poisson}.

\begin{proposition} \label{poissongen} Let the conditions and notations of Proposition \ref{poisson} be kept and assume that
$B>0$ and $a\in \mathbb{R}^d$. Then
$$
\sum\limits_{y\in a+\Lambda} f\left(\frac{y}{B}\right) = \frac{B^d}{\mbox{\rm Vol}(\mathbb{R}^d/\Lambda)} \cdot 
\sum\limits_{x\in \Lambda'} e(a\cdot x)\hat{f}(Bx).
$$
\end{proposition} 
$ $\\
{\bf Proof of Theorem \ref{ls}:} We first note that
\begin{equation} \label{KK}
\begin{split}
K = & \max\limits_{1\le i\le R} \sharp \left\{j\in \{1,...,R\} : \min\limits_{z\in \mathbb{Z}^d} ||x_j-x_i-z||_2\le \sqrt{d}N^{-1/d}\right\} \\
  \geq & \max\limits_{1\le i\le R} \sharp \{j\in \{1,...,R\} : \min\limits_{z\in \mathbb{Z}^d} \max\limits_{1\le k\le d} |x_j^{(k)}-x_i^{(k)}-z^{(k)}| 
  \le N^{-1/d}\}\\
  = & \max\limits_{1\le i\le R} \sharp \{j\in \{1,...,R\} : \max\limits_{1\le k\le d} ||x_j^{(k)}-x_i^{(k)}-z^{(k)}|| \le N^{-1/d}\} =: K',
\end{split}
\end{equation}
where $||u||$ is the distance of $u\in \mathbb{R}$ to the nearest integer and we write
$$
x_i=\left(x_i^{(1)},...,x_i^{(d)}\right)\quad \mbox{and} \quad z=\left(z^{(1)},...,z^{(d)}\right)
$$
for $i=1,...,R$. 

Now let  $S = \{x_1, x_2, . . . , x_R\}$. Taking Proposition \ref{duality}, the duality principle, into account, it suffices to prove that 
$$
\sum\limits_{\substack{n\in \mathbb{Z}^d\\ ||n||_2\le N^{1/d}}}\left|\sum\limits_{x\in S} b_x\cdot e(n\cdot x)\right|^2 \ll KN\sum\limits_{x\in S}|b_x|^2
$$ 
for any complex numbers $b_x$. To this end, for $x=\left(x^{(1)},...,x^{(d)}\right)\in \mathbb{R}^d$, we define 
\begin{equation*}
      \phi(x) = \prod\limits_{k=1}^d \left(\frac{\sin\left(\pi x^{(k)}\right)}{2x^{(k)}}\right)^2
\end{equation*}
and note that $\phi(x)$ is non-negative and satisfies $\phi(x)\ge 1$ if $|x^{(k)}|\le 1/2$ for $k=1,...d$. 
Moreover, the Fourier transformation of $\phi(x)$ equals
$$
\hat{\phi}(s) = \left(\frac{\pi^2}{4}\right)^d\prod\limits_{k=1}^d\max\left\{1-\left|s^{(k)}\right|,\ 0\right\},
$$
where we write 
$$
s=\left(s^{(1)},...,s^{(d)}\right).
$$
Hence,
\begin{equation*}
\begin{split}
\sum\limits_{\substack{n\in \mathbb{Z}^d\\ ||n||_2\le N^{1/d}}} 
\left|\sum\limits_{x\in S} b_x\cdot e(n\cdot x)\right|^2 \le & \sum\limits_{n\in\mathbb{Z}^d}\phi\left(\frac{n}{2N^{1/d}}\right)
                                                       \left|\sum\limits_{x\in S} b_x\cdot e(n\cdot x)\right|^2\\
                             = & \sum\limits_{x,x'\in S} b_x \overline{b_{x'}}\cdot V(x-x'),
\end{split}
\end{equation*}
where
$$
V(y) = \sum\limits_{n\in \mathbb{Z}^d} \phi\left(\frac{n}{2N^{1/d}}\right)\cdot e(n\cdot y).
$$
Using Proposition \ref{poissongen}, the Poisson summation formula, we transform $V(y)$ into
\begin{equation*}
\begin{split}
V(y) = & 2^d N \cdot \sum\limits_{\alpha \in y+\mathbb{Z}^d} \tilde{\phi}\left(2N^{1/d} \alpha\right)
     = 2^d N \cdot \sum\limits_{\alpha\in -y+\mathbb{Z}^d} \hat{\phi}\left(2N^{1/d} \alpha \right)\\
     = & \frac{\pi^{2d}}{2^d}\cdot N \cdot \prod\limits_{k=1}^d \max\left\{1-\left|2N^{1/d}y^{(k)}\right|,0\right\},
\end{split}
\end{equation*}
where $\tilde{\phi}$ is the inverse Fourier transform and $\hat{\phi}$ is the Fourier transform of $\phi$.
Therefore,
\begin{equation}
\begin{split}
\sum\limits_{\substack{n\in \mathbb{Z}^d\\ ||n||_2\le N^{1/d}}}\left|\sum\limits_{x\in S} b_x\cdot e(n\cdot x)\right|^2 \le 
\frac{\pi^{2d}}{2^d}\cdot N \cdot \sum\limits_{\substack{x,x'\in S\\ 
||x^{(i)}-{x'}^{(i)}||\le N^{-1/d} \mbox{\scriptsize\ for\ } i=1,...,d}} |b_x| |b_{x'}|.
\end{split}
\end{equation}
Now we observe that
$$
|b_x| |b_{x'}| \le \frac{1}{2}\cdot \left(|b_x|^2+|b_{x'}|^2\right). 
$$
It follows that 
\begin{equation*}
\begin{split}
\sum\limits_{\substack{n\in \mathbb{Z}^d\\ ||n||_2\le N^{1/d}}}\left|\sum\limits_{x\in S} b_x\cdot e(n\cdot x)\right|^2 \ll 
K'NZ\le KNZ,
\end{split}
\end{equation*}
where we use \eqref{Zdef} and \eqref{KK}.
This completes the proof. $\Box$

\section{Conversion into a counting problem}
Now we return to the large sieve for $\mathbb{Q}(i)$. 
We begin with restricting the moduli $q$ to an arbitrary multiset $S$ of elements of $\mathbb{Z}[i]\setminus \{0\}$. We shall also restrict 
the norms of these moduli to dyadic intervals, which is for technical reasons. Thus, we are interested in 
estimating the quantity
$$
T:=\sum\limits_{\substack{q\in S\\ Q/2<\mathcal{N}(q)\le Q}} 
\sum\limits_{\substack{r \bmod{q}\\ (r,q)=1}} \left|\sum\limits_{\substack{n\in \mathbb{Z}[i]\\ \mathcal{N}(n)\le N}}  
a_n \cdot e\left(\Re\left(\frac{nr}{q}\right)\right)\right|^2. 
$$
We shall later confine ourselves to squares or, more generally, $k$-th powers. 

Our first step is to re-write $T$ in the form
\begin{equation} \label{T}
T=\sum\limits_{\substack{q\in S\\ Q/2<\mathcal{N}(q)\le Q}} 
\sum\limits_{\substack{r \bmod{q}\\ (r,q)=1}} \left|\sum\limits_{\substack{n\in \mathbb{Z}[i]\\ \mathcal{N}(n)\le N}}  
a_n \cdot e\left(\left(\frac{xu+yv}{\mathcal{N}(q)},\frac{xv-yu}{\mathcal{N}(q)}\right)\cdot (s,t)\right)\right|^2, 
\end{equation}
where
$$
q=u+iv,\quad r=x+iy,\quad n=s+ti. 
$$
To bound $T$, we employ Theorem \ref{ls} for the case $d=2$, which immediately gives us the following.

\begin{corollary} \label{preform} For $T$ as defined in \eqref{T}, we have the bound
$$
T\ll KNZ,
$$
where $Z$ is defined as in \eqref{Zdef} and 
\begin{equation*}
\begin{split}
& K:=\max\limits_{r_1,q_1} \sharp \Bigg\{(r_2,q_2) : \\
& \min\limits_{z\in \mathbb{Z}^2} 
\left|\left|\left(\frac{x_2u_2+y_2v_2}{\mathcal{N}(q_2)},
\frac{x_2v_2-y_2u_2}{\mathcal{N}(q_2)}\right)-
\left(\frac{x_1u_1+y_1v_1}{\mathcal{N}(q_1)},\frac{x_1v_1-y_1u_1}{\mathcal{N}(q_1)}\right)-z\right|\right|_2^2
\le 2N^{-1}\Bigg\}
\end{split}
\end{equation*}
with the conventions that, for $j=1,2$, 
$$
q_j\in S,\quad Q/2<\mathcal{N}(q_j)\le Q,
$$
$\{r_j\}$ forms a system of representatives of reduced residue classes modulo $q_j$
and 
$$
q_j=u_j+iv_j,\quad r_j=x_j+iy_j. 
$$
\end{corollary}

Thus, we have converted the problem into a counting problem.

\section{Switching back to $\mathbb{Z}[i]$}
Now we observe that
$$
\frac{\overline{r_j}}{\overline{q_j}}=\frac{x_j-iy_j}{u_j-iv_j}=\frac{x_ju_j+y_jv_j}{\mathcal{N}(q_j)}+\frac{x_jv_j-y_ju_j}{\mathcal{N}(q_j)}i.
$$
It follows that
\begin{equation*}
\begin{split}
K= & \max\limits_{r_1,q_1} \sharp \Bigg\{(r_2,q_2) : 
\min\limits_{z\in \mathbb{Z}[i]} 
\left|\frac{\overline{r_2}}{\overline{q_2}}-\frac{\overline{r_1}}{\overline{q_1}}-z \right|^2
\le 2N^{-1}\Bigg\}\\
= & \max\limits_{\substack{q_1\in S\\ Q/2<\mathcal{N}(q_1)\le Q\\ (r_1,q_1)=1}} 
\sharp \Bigg\{(r_2,q_2) : q_2\in S,\ Q/2<\mathcal{N}(q_2)\le Q,\ (r_2,q_2)=1, \\ & 
\left|\frac{r_2}{q_2}-\frac{r_1}{q_1}\right|^2
\le 2N^{-1}\Bigg\}.
\end{split}
\end{equation*}
Further,
$$
\left|\frac{r_2}{q_2}-\frac{r_1}{q_1}\right|_2^2
\le 2N^{-1} \Longleftrightarrow \mathcal{N}\left(r_1q_2-r_2q_1\right)
\le 2N^{-1}\mathcal{N}(q_1)\mathcal{N}(q_2). 
$$
We deduce that 
\begin{equation} \label{Ktransform}
\begin{split}
K\le & \max\limits_{\substack{q_1\in S\\ Q/2<\mathcal{N}(q_1)\le Q\\ (r_1,q_1)=1}} \ 
\sum\limits_{\substack{b\in \mathbb{Z}[i]\\ \mathcal{N}(b)\le 2N^{-1}\mathcal{N}(q_1)\mathcal{N}(q_2)}} \
\sum\limits_{\substack{q_2\in S,\\ Q/2<\mathcal{N}(q_2)\le Q\\ r_1q_2\equiv b \bmod{q_1}}} 1 \\
\ll & \max\limits_{\substack{q_1\in S\\ Q/2<\mathcal{N}(q_1)\le Q\\ (r_1,q_1)=1}} \ 
\sum\limits_{b\in \mathbb{Z}[i]} \Phi_1\left(\mathcal{N}\left(\frac{b\sqrt{N}}{q_1\sqrt{2Q}}\right)\right) 
\sum\limits_{\substack{q_2\in S,\\ r_1q_2\equiv b \bmod{q_1}}} \Phi_2\left(\mathcal{N}\left(\frac{q_2}{\sqrt{Q}}\right)\right)\\
= & \max\limits_{\substack{q_1\in S\\ Q/2<\mathcal{N}(q_1)\le Q\\ (r_1,q_1)=1}} \
\sum\limits_{q_2\in S} \Phi_2\left(\mathcal{N}\left(\frac{q_2}{\sqrt{Q}}\right)\right)\cdot 
\sum\limits_{b\equiv r_1q_2 \bmod{q_1}} \Phi_1\left(\mathcal{N}\left(\frac{b\sqrt{N}}{q_1\sqrt{2Q}}\right)\right),
\end{split}
\end{equation}
where, for $i=1,2$, $\Phi_{i} : \mathbb{R} \rightarrow \mathbb{R}^+$ are any smooth functions with rapid decay such that 
$\Phi_{i}(x)\gg 1$ if $|x|\le 1$. 
We shall fix $\Phi_{i}$ later suitably. Let $\Psi_{i}: \mathbb{C}\rightarrow \mathbb{R}^+$ be given by 
$$
\Psi=\Phi \circ \mathcal{N} \quad \mbox{for } i=1,2.
$$ 
Then the above inequality for $K$ turns into
\begin{equation} \label{superK}
\begin{split}
K\ll \max\limits_{\substack{q_1\in S\\ Q/2<\mathcal{N}(q_1)\le Q\\ (r_1,q_1)=1}} \
\sum\limits_{q_2\in S} \Psi_2\left(\frac{q_2}{\sqrt{Q}}\right)\cdot
\sum\limits_{b\equiv r_1q_2 \bmod{q_1}} \Psi_1\left(\frac{b\sqrt{N}}{|q_1|\sqrt{2Q}}\right).
\end{split}
\end{equation}

\section{Application of Poisson summation}
To transform the inner-most sum over $b$, we use Poisson summation again. The complex numbers $a\equiv 0 \bmod{q_1}$ form a square lattice
$$
\Lambda=\left\{x\binom{u_1}{v_1}+y\binom{-v_1}{u_1}: (x,y)\in \mathbb{Z}^2\right\}\subset \mathbb{R}^2
$$
with volume $\mathcal{N}(q_1)$ when regarded as vectors in $\mathbb{R}^2$. The dual lattice turns out to be 
$$
\Lambda'=\frac{1}{\mathcal{N}(q_1)}\cdot \Lambda,
$$
which corresponds to the set 
$$
\left\{\frac{a}{\mathcal{N}(q_1)} : a\equiv 0 \bmod{q_1}\right\}
$$
in $\mathbb{C}$. Hence, Proposition \ref{poisson} gives
\begin{equation} \label{afterpoisson}
\begin{split}
\sum\limits_{b\equiv r_1q_2 \bmod{q_1}} \Psi_1\left(\frac{b\sqrt{N}}{|q_1|\sqrt{2Q}}\right)= & 
\frac{2Q}{N}\cdot \sum\limits_{a\equiv 0\bmod{q_1}} e\left( \frac{\overrightarrow{a}\cdot \overrightarrow{r_1q_2}}{\mathcal{N}(q_1)}\right) 
\hat\Psi_1\left(\frac{a|q_1|\sqrt{2Q}}{\mathcal{N}(q_1)\sqrt{N}}\right)\\
= & \frac{2Q}{N}\cdot \sum\limits_{j\in \mathbb{Z}[i]} e\left(\overrightarrow{\frac{j}{\overline{q_1}}}\cdot \overrightarrow{r_1q_2}\right) 
\hat\Psi_1\left(\frac{j\sqrt{2Q}}{\sqrt{N}}\right),
\end{split}
\end{equation}
where for $x,z\in \mathbb{C}$, we write
$$
\overrightarrow{z}=\binom{\Re(z)}{\Im(z)} \quad \mbox{and} \quad \hat\Psi(x)=\int\limits_{\mathbb{C}} \Psi(y) e\left(-\overrightarrow{y}\cdot
\overrightarrow{x}\right) \ dy_2dy_1
$$
with $y_1:=\Re(y)$ and $y_2:=\Im(y)$. 
Combining \eqref{superK} and \eqref{afterpoisson}, and re-arranging summation, we deduce that
\begin{equation} \label{Kgen}
K\ll \frac{Q}{N}\cdot\max\limits_{\substack{q_1\in S\\ Q/2<\mathcal{N}(q_1)\le Q\\ (r_1,q_1)=1}} \
\sum\limits_{j\in \mathbb{Z}[i]} \hat\Psi_1\left(\frac{j\sqrt{2Q}}{\sqrt{N}}\right) \cdot
\sum\limits_{q_2\in S} \Psi_2\left(\frac{q_2}{\sqrt{Q}}\right)\cdot 
e\left(\overrightarrow{\frac{j}{\overline{q_1}}}\cdot \overrightarrow{r_1q_2}\right).
\end{equation}
We observe that for $a,b\in \mathbb{C}$,
\begin{equation} \label{obs}
e\left(\overrightarrow{a}\cdot \overrightarrow{b}\right)=e\left(\Re(\overline{a}b)\right). 
\end{equation}
Hence, upon a change of variables $j \rightarrow \overline{j}$, we arrive at
\begin{equation} \label{Kgen1}
K\ll \frac{Q}{N}\cdot\max\limits_{\substack{q_1\in S\\ Q/2<\mathcal{N}(q_1)\le Q\\ (r_1,q_1)=1}} \
\sum\limits_{j\in \mathbb{Z}[i]} \hat\Psi_1\left(\frac{j\sqrt{2Q}}{\sqrt{N}}\right) 
\sum\limits_{q_2\in S} \Psi_2\left(\frac{q_2}{\sqrt{Q}}\right)\cdot
e\left(\Re\left(\frac{jr_1}{q_1}\cdot q_2\right)\right).
\end{equation}

\section{The case of square moduli}
Now we restrict overselves to the case when $S$ is the set of non-zero squares in $\mathbb{Z}[i]$. We write $Q_0=\sqrt{Q}$
and replace $q_i$ by $q_i^2$ ($i=1,2$). Throughout the following, we assume that $Q_0>N^{1/4}$ for otherwise the desired result follows from 
\eqref{Huxley} upon extending the set of moduli to all non-zero Gaussian integers. 
We deduce from \eqref{Kgen1} that 
\begin{equation} \label{Knew}
\begin{split}
K\ll & \frac{Q_0^2}{N}\cdot \max\limits_{\substack{Q_0/\sqrt{2}<\mathcal{N}(q_1)\le Q_0\\ (r_1,q_1)=1}} 
\sum\limits_{j\in \mathbb{Z}[i]} \hat\Psi_1\left(\frac{\sqrt{2}jQ_0}{\sqrt{N}}\right) \cdot  
\sum\limits_{q_2\in \mathbb{Z}[i]} \Psi_2\left(\frac{q_2^2}{Q_0}\right)\cdot 
e\left(\Re\left(\frac{jr_1}{q_1^2}\cdot q_2^2 \right) \right)\\
\ll & \frac{Q_0^3}{N}+ \frac{Q_0^2}{N}\cdot \max\limits_{\substack{Q_0/\sqrt{2}<\mathcal{N}(q_1)\le Q_0\\ (r_1,q_1)=1}} 
\sum\limits_{j\in \mathbb{Z}[i]\setminus \{0\}} \hat\Psi_1\left(\frac{\sqrt{2}jQ_0}{\sqrt{N}}\right) S\left(q_1,r_1,j\right),
\end{split}
\end{equation}
where 
$$
S\left(q_1,r_1,j\right):=\sum\limits_{q_2\in \mathbb{Z}[i]} \Psi_2\left(\frac{q_2^2}{Q_0}\right) \cdot 
e\left(\Re\left(\frac{jr_1}{q_1^2}\cdot q_2^2\right)\right).
$$
Here we use the estimate 
$$
\sum\limits_{q_2\in \mathbb{Z}[i]} \Psi_2\left(\frac{q_2^2}{Q_0}\right) \ll Q_0
$$
to bound the contribution of $j=0$.  

\subsection{Weyl differencing}
Now we employ Weyl differencing in the setting of $\mathbb{Z}[i]$. 
Using the Cauchy-Schwarz inequality, we deduce that
\begin{equation} \label{afterCS}
\begin{split}
K\ll & \frac{Q_0^3}{N}+ \frac{Q_0^2}{N}\cdot \max\limits_{\substack{Q_0/\sqrt{2}<\mathcal{N}(q_1)\le Q_0\\ (r_1,q_1)=1}} 
\left(\sum\limits_{j\in \mathbb{Z}[i]} \hat\Psi_1\left(\frac{\sqrt{2}jQ_0}{\sqrt{N}}\right)\right)^{1/2}\times\\ 
& \left(\sum\limits_{j\in \mathbb{Z}[i]\setminus \{0\}} \hat\Psi_1\left(\frac{jQ_0}{\sqrt{N}}\right)\cdot \left|
S\left(q_1,r_1,j\right) \right|^2\right)^{1/2}\\
\ll & \frac{Q_0^3}{N}+ \frac{Q_0}{\sqrt{N}}\cdot \max\limits_{\substack{Q_0/\sqrt{2}<\mathcal{N}(q_1)\le Q_0\\ (r_1,q_1)=1}}  
\left(\sum\limits_{j\in \mathbb{Z}[i]\setminus \{0\}} \hat\Psi_1\left(\frac{\sqrt{2}jQ_0}{\sqrt{N}}\right)\cdot \left|
S\left(q_1,r_1,j\right) \right|^2\right)^{1/2},
\end{split}
\end{equation}
where we use the estimate
$$
\sum\limits_{j\in \mathbb{Z}[i]} \hat\Psi_1\left(\frac{\sqrt{2}jQ_0}{\sqrt{N}}\right) \ll \frac{N}{Q_0^2}.
$$

Multiplying out the square, we get
\begin{equation} \label{multi}
\left| S\left(q_1,r_1,j\right) \right|^2 =  
\sum\limits_{q_2,q\in \mathbb{Z}[i]} \Psi_2\left(\frac{q_2^2}{Q_0}\right) \cdot \Psi_2\left(\frac{q^2}{Q_0}\right) \cdot   
e\left(\Re\left(\frac{jr_1}{q_1^2}\cdot \left(q_2^2-q^2\right) \right)\right).
\end{equation}
We now set
$$
\alpha:=q_2-q
$$
so that 
$$
q_2^2-q^2=\alpha^2+2\alpha q \quad \mbox{and} \quad q_2=\alpha+q.
$$
Then \eqref{multi} turns into
\begin{equation} \label{change}
\begin{split}
\left| S\left(q_1,r_1,j\right) \right|^2 = &
\sum\limits_{\alpha\in \mathbb{Z}[i]} e\left(\Re\left(\frac{jr_1\alpha^2}{q_1^2}\right)\right) \times\\ &
\sum\limits_{q\in \mathbb{Z}[i]} \Psi_2\left(\frac{q^2}{Q_0}\right) \cdot\Psi_2\left(\frac{(\alpha+q)^2}{Q_0}\right) \cdot  
e\left(\Re\left(\frac{2j\alpha r_1}{q_1^2}\cdot q\right)\right).
\end{split}
\end{equation}

\subsection{Poisson summation} 
We shall apply Proposition \ref{poissongen} with $d=2$ to transform the inner-most over $q$ on the right-hand side of \eqref{change}. For 
$z=(z_1,z_2)\in \mathbb{R}^2$, we set 
$$
g(z):=\Psi_2\left((z_1+iz_2)^2\right) \cdot\Psi_2\left(\left(\frac{\alpha}{\sqrt{Q_0}}+z_1+iz_2\right)^2\right).
$$
Then using \eqref{obs} with 
$$
a=q \quad \mbox{and} \quad b=\overline{\frac{2j\alpha r_1}{q_1^2}},
$$
we deduce that
\begin{equation} \label{backtor2}
\sum\limits_{q\in \mathbb{Z}[i]} \Psi_2\left(\frac{q^2}{Q_0}\right) \cdot\Psi_2\left(\frac{(\alpha+q)^2}{Q_0}\right) \cdot  
e\left(\Re\left(\frac{2j\alpha r_1}{q_1^2}\cdot q\right)\right)
=\sum\limits_{x\in \mathbb{Z}^2} e\left(\overrightarrow{b}\cdot x\right)g\left(\frac{x}{\sqrt{Q_0}}\right).
\end{equation}
Now applying Proposition \ref{poissongen} with 
$B:=1/\sqrt{Q_0}$ and $f:=\tilde{g}$, the inverse Fourier transform of $g$, to the right-hand side of \eqref{backtor2}, we get
\begin{equation*}
\sum\limits_{x\in \mathbb{Z}^2} e\left(\overrightarrow{b}\cdot x\right)g\left(\frac{x}{\sqrt{Q_0}}\right) 
= Q_0\cdot \sum\limits_{y\in \overrightarrow{b}+\mathbb{Z}^2} \tilde{g}\left(\sqrt{Q_0}y\right) =
Q_0\cdot \sum\limits_{y\in -\overrightarrow{b}+\mathbb{Z}^2} \hat{g}\left(\sqrt{Q_0}y\right).
\end{equation*}
It follows that
\begin{equation} \label{newpoiss}
\left| S\left(q_1,r_1,j\right) \right|^2
= Q_0\cdot \sum\limits_{\alpha\in \mathbb{Z}[i]} e\left(\Re\left(\frac{jr_1\alpha^2}{q_1^2}\right)\right) \cdot 
\sum\limits_{\beta\in \mathbb{Z}[i]} \hat{g}\left(\sqrt{Q_0}\cdot 
\overrightarrow{\left(\beta-\overline{\frac{2j\alpha r_1}{q_1^2}}\right)}\right).
\end{equation}

At this point, we specify our choice of $\Psi_2$ and compute the Fourier transform of $g$. We set 
$$
\Phi_2(t):=\exp\left(-\frac{\pi}{2}\cdot \sqrt{|t|}\right)
$$
so that 
$$
\Psi_2(z)=\Phi_2(\mathcal{N}(z))=\exp\left(-\frac{\pi}{2}\cdot \sqrt{\mathcal{N}(z)}\right).
$$
It follows that
$$
g(z)=\exp\left(-\frac{\pi}{2}\cdot \left(z_1^2+z_2^2+\left(\frac{\alpha_1}{\sqrt{Q_0}}+z_1\right)^2+
\left(\frac{\alpha_2}{\sqrt{Q_0}}+z_2\right)^2\right)\right),
$$
where $\alpha_1:=\Re(\alpha)$ and $\alpha_2:=\Im(\alpha)$. Completing the squares, it follows that
$$
g(z)=\exp\left(-\frac{\pi}{4Q_0}\cdot \left(\alpha_1^2+\alpha_2^2\right)\right)\cdot
\exp\left(-\pi\left(\left(z_1+\frac{\alpha_1}{2\sqrt{Q_0}}\right)^2+\left(z_2+\frac{\alpha_2}{2\sqrt{Q_0}}\right)^2\right)\right).
$$
The Fourier transform of this function equals
\begin{equation} \label{gfourier}
\begin{split}
\hat{g}(z)= & \exp\left(-\frac{\pi}{4Q_0}\cdot \left(\alpha_1^2+\alpha_2^2\right)\right)\cdot e\left(-\frac{\alpha_1z_1+\alpha_2z_2}{2\sqrt{Q_0}}\right)\cdot
\exp\left(-\pi \left(z_1^2+z_2^2\right)\right)\\ = &\exp\left(-\frac{\pi}{4Q_0}\cdot \mathcal{N}(\alpha)\right)\cdot 
e\left(-\frac{\Re(\overline{\alpha}(z_1+iz_2))}{2\sqrt{Q_0}}\right)\cdot
\exp\left(-\pi \mathcal{N}(z_1+iz_2)\right).
\end{split}
\end{equation}
Plugging \eqref{gfourier} into \eqref{newpoiss}, using the triangle inequality and bounding all terms of the form $e(\gamma)$ trivially by 
$|e(\gamma)|\le 1$, we get
\begin{equation} \label{finalafterps}
\begin{split}
\left| S\left(q_1,r_1,j\right) \right|^2\le & Q_0 \cdot \sum\limits_{\alpha\in \mathbb{Z}[i]} \exp\left(-\frac{\pi}{4Q_0}\cdot \mathcal{N}(\alpha)\right)
\times\\ &
\sum\limits_{\beta\in \mathbb{Z}[i]} \exp\left(-\pi Q_0 \mathcal{N}\left(\beta-\overline{\frac{2j \alpha r_1}{q_1^2}}\right)\right).         
\end{split}
\end{equation}

\subsection{Counting}
The contributions of $\beta$'s with 
$$
\mathcal{N}\left(\beta-\overline{\frac{2j \alpha r_1}{q_1^2}}\right) > Q_0^{\varepsilon-1}
$$
and of $\alpha$'s with 
$$
\mathcal{N}(\alpha)> Q_0^{1+\varepsilon}
$$
to the right-hand side of \eqref{finalafterps} are neglible. Therefore, it follows from \eqref{finalafterps} that
\begin{equation} \label{squarebound}
 \begin{split}
     \left| S\left(q_1,r_1,j\right) \right|^2
        \ll & \frac{1}{(Q_0N)^{2018}}+Q_0\sum\limits_{\mathcal{N}(\alpha)\le Q_0^{1+\varepsilon}} \sum\limits_{\substack{\beta\in \mathbb{Z}[i]\\
        \mathcal{N}\left(\beta-\overline{2j \alpha r_1/q_1^2}\right)\le Q_0^{\varepsilon-1}}} 1\\
        \ll & \frac{1}{(Q_0N)^{2018}}+Q_0\sum\limits_{\substack{\mathcal{N}(\alpha)\le Q_0^{1+\varepsilon}\\ ||2j\alpha r_1/q_1^2||\le Q_0^{\varepsilon-1/2}}} 1,        
 \end{split}
\end{equation}
where $||z||$ is the distance of $z\in \mathbb{C}$ to the nearest Gaussian integer. 

Now we want to bound the term in the maximum on the right-hand side of \eqref{afterCS}. To this end, we choose $\Psi_1$ 
in a suitable way so that $\hat\Psi_1$ decays exponentially. We set 
$$
\Phi_1(t):=\exp\left(-\pi |t|\right)
$$
Since $\Psi_1=\Phi_1\circ \mathcal{N}$, it follows that
\begin{equation} \label{psi1}
\Psi_1(z)=\exp\left(-\pi \mathcal{N}(z)\right)=\hat{\Psi}_1(z).
\end{equation}
Hence, using \eqref{squarebound}, we obtain
\begin{equation} \label{obvious}
\sum\limits_{j\in \mathbb{Z}[i]\setminus \{0\}} \hat\Psi_1\left(\frac{\sqrt{2}jQ_0}{\sqrt{N}}\right)\cdot \left| S\left(q_1,r_1,j\right) \right|^2\\
\ll 1+Q_0\sum\limits_{0<\mathcal{N}(j)\le NQ_0^{\varepsilon-2}} 
\sum\limits_{\substack{\mathcal{N}(\alpha)\le Q_0^{1+\varepsilon}\\ ||2j\alpha r_1/q_1^2||\le Q_0^{\varepsilon-1/2}}} 1
\end{equation}
upon noting that the contribution of $\mathcal{N}(j)>NQ_0^{\varepsilon-2}$ is negligible. 
The contribution of $\alpha=0$ to the right-hand side of \eqref{obvious} is obviously bounded by
$$
\ll \frac{N}{Q_0^{1-\varepsilon}}.
$$
Writing $d=2j\alpha$ and noting that the number of divisors of $d\in \mathbb{Z}[i]\setminus \{0\}$ in the Gaussian integers is bounded by 
$O\left(\mathcal{N}(d)^{\varepsilon}\right)$, we deduce that
\begin{equation} \label{super}
\begin{split}
\sum\limits_{j\in \mathbb{Z}[i]\setminus \{0\}} \hat\Psi_1\left(\frac{\sqrt{2}jQ_0}{\sqrt{N}}\right)\cdot \left| S\left(q_1,r_1,j\right) \right|^2
\ll & \frac{N}{Q_0^{1-\varepsilon}}+ 
N^{\varepsilon}Q_0\sum\limits_{\substack{\mathcal{N}(d)\le 4NQ_0^{2\varepsilon-1}\\ ||d r_1/q_1^2||\le Q_0^{\varepsilon-1/2}}} 1\\
\ll & \frac{N}{Q_0^{1-\varepsilon}}+ 
N^{\varepsilon}Q_0\sum\limits_{\substack{l\in \mathbb{Z}[i]\\ |l/q_1^2|\le Q_0^{\varepsilon-1/2}}} 
\sum\limits_{\substack{\mathcal{N}(d)\le 4NQ_0^{2\varepsilon-1}\\ d\equiv l\overline{r_1} \bmod{q_1^2}}} 1, 
\end{split}
\end{equation}
where $\overline{r_1}$ is a multiplicative inverse of $r_1$ modulo $q_1^2$, i.e. $r_1\overline{r_1}\equiv 1 \bmod{q_1^2}$. 
The number of residue classes modulo $q_1^{2}$ is $\mathcal{N}(q_1^2)\le Q_0^2$, and hence
\begin{equation} \label{res}
\sum\limits_{\substack{\mathcal{N}(d)\le 4NQ_0^{2\varepsilon-1}\\ d\equiv l\overline{r_1} \bmod{q_1^2}}} 1 \ll 1+\frac{N}{Q_0^{3-2\varepsilon}}.
\end{equation}
Further,
\begin{equation} \label{lat}
\sum\limits_{\substack{l\in \mathbb{Z}[i]\\ 
|l/q_1^2|\le Q_0^{\varepsilon-1/2}}} 1 \le \sum\limits_{\substack{l\in \mathbb{Z}[i]\\ 
|l|\le Q_0^{\varepsilon+1/2}}} 1 \ll Q_0^{1+2\varepsilon}.
\end{equation}
Combining \eqref{super}, \eqref{res} and \eqref{lat}, we obtain
\begin{equation} \label{comb1}
\sum\limits_{j\in \mathbb{Z}[i]\setminus \{0\}} \hat\Psi_1\left(\frac{\sqrt{2}jQ_0}{\sqrt{N}}\right)\cdot \left|
S\left(q_1,r_1,j\right)\right|^2
\ll (Q_0N)^{4\varepsilon}\left( Q_0^{2}+\frac{N}{Q_0} \right),
\end{equation}
and combining \eqref{afterCS} and \eqref{comb1}, we arrive at
\begin{equation} \label{comb2}
K\ll \frac{Q_0^3}{N}+(Q_0N)^{2\varepsilon}\left( \frac{Q_0^2}{N^{1/2}}+Q_0^{1/2} \right),
\end{equation}
our final bound for $K$. 

Now the statement in Theorem \ref{squaremodZi} follows immediately from Corollary \ref{preform} and 
\eqref{comb2} after dividing the moduli into dyadic intervals and replacing $Q_0$ by $Q$.

\section{The case of $k$-th power moduli}
Now we take $S$ as the set of non-zero $k$th-powers in $\mathbb{Z}[i]$. We write $Q_0=Q^{1/k}$
and replace $q_i$ by $q_i^k$ ($i=1,2$).  Throughout the following, we assume that $Q_0>N^{1/(2k)}$ for otherwise the desired result follows from 
\eqref{Huxley} upon extending the set of moduli to all non-zero Gaussian integers.
We deduce from \eqref{Kgen1} that
\begin{equation} \label{generalKbound} 
\begin{split}
K\ll & \frac{Q_0^k}{N}\cdot \max\limits_{\substack{Q_0/\sqrt[k]{2}<\mathcal{N}(q_1)\le Q_0\\ (r_1,q_1)=1}} 
\sum\limits_{j\in \mathbb{Z}[i]} \hat\Psi_1\left(\frac{\sqrt{2}jQ_0^{k/2}}{\sqrt{N}}\right) \times\\ &  
\sum\limits_{q_2\in \mathbb{Z}[i]} \Psi_2\left(\frac{q_2^k}{Q_0^{k/2}}\right)\cdot 
e\left(\Re\left(\frac{jr_1}{q_1^k}\cdot q_2^k \right) \right)\\
\ll & \frac{Q_0^{k+1}}{N}+ \frac{Q_0^k}{N}\cdot \max\limits_{\substack{Q_0/\sqrt[k]{2}<\mathcal{N}(q_1)\le Q_0\\ (r_1,q_1)=1}} 
\sum\limits_{j\in \mathbb{Z}[i]\textbackslash \{0\}} \hat\Psi_1\left(\frac{\sqrt{2}jQ_0^{k/2}}{\sqrt{N}}\right) \cdot 
\left| S_k\left(q_1,r_1,j\right)\right|,
\end{split}
\end{equation}
where 
$$
S_k\left(q_1,r_1,j\right)=\sum\limits_{q_2\in \mathbb{Z}[i]} \Psi_2\left(\frac{q_2^k}{Q_0^{k/2}}\right) \cdot 
e\left(\Re\left(\frac{jr_1}{q_1^k}\cdot q_2^k \right) \right).
$$
Here we use the estimate 
$$
\sum\limits_{q_2\in \mathbb{Z}[i]} \Psi_2\left(\frac{q_2^k}{Q_0^{k/2}}\right) \ll Q_0
$$
to bound the contribution of $j=0$. 

\subsection{Weyl differencing}
Multiplying out the square and setting $\alpha_1=q_2-q$, we obtain
\begin{equation*}
\begin{split}
& \left|S_k\left(q_1,r_1,j\right) \right|^2 \\
=  & \sum\limits_{q_2,q\in \mathbb{Z}[i]} \Psi_2\left(\frac{q_2^k}{Q_0^{k/2}}\right) \cdot \Psi_2\left(\frac{q^k}{Q_0^{k/2}}\right) \cdot 
e\left(\Re\left(\frac{jr_1}{q_1^k}\cdot(q_2^k-q^k) \right) \right),\\
 =  & \sum\limits_{\alpha_1,q\in \mathbb{Z}[i]} \Psi_2\left(\frac{q^k}{Q_0^{k/2}}\right) \cdot \Psi_2\left(\frac{(\alpha_1+q)^k}{Q_0^{k/2}}\right) 
 \cdot e\left(\Re\left(\frac{jr_1}{q_1^k}\cdot\left((\alpha_1 + q)^k-q^k\right) \right) \right).
 \end{split}
 \end{equation*}
We observe that the contribution of $\alpha_1$'s with $\mathcal{N}\left(\alpha_1\right)>Q_0^{1+\varepsilon}$ is negligible and 
write
$$
P_{k-1,\alpha_1}(q) = (\alpha_1 + q)^k - q^k = \binom{k}{1} \cdot \alpha_1 q^{k-1} + \binom{k}{2}\cdot \alpha_1^2q^{k-2} +\cdots +
\binom{k}{k}\cdot \alpha_1^k,
$$
thus obtaining
\begin{equation*}
\left|S_k\left(q_1,r_1,j\right) \right|^2 
 \ll \left|\sum\limits_{\mathcal{N}(\alpha_1)\le Q_0^{1+\varepsilon}} S_{k-1}\left(q_1,r_1,j,\alpha_1\right)\right|,
 \end{equation*}
 where
 $$
 S_{k-1}\left(q_1,r_1,j,\alpha_1\right):=\sum\limits_{q\in \mathbb{Z}[i]} 
\Psi_2\left(\frac{q^k}{Q_0^{k/2}}\right) \cdot \Psi_2\left(\frac{(\alpha_1+q)^k}{Q_0^{k/2}}\right)  
\cdot e\left(\Re\left(\frac{jr_1}{q_1^k}\cdot P_{k-1,\alpha_1}(q)\right) \right).
 $$
 If $k > 2$, we apply the Cauchy-Schwarz inequality again to obtain
 \begin{equation*}
 \left|S_k\left(q_1,r_1,j\right) \right|^4 \ll
 Q_0^{1+ \varepsilon}\sum\limits_{\substack{\alpha_1 \in \mathbb{Z}[i] \\ \mathcal{N}(\alpha_1) \le Q_0^{1+\varepsilon}}} 
 \left| S_{k-1}\left(q_1,r_1,j,\alpha_1\right)
 \right|^2. 
 \end{equation*}
Multiplying out the square, changing variables and truncating the resulting sums in a similar way as above, we now obtain
\begin{equation*}
\begin{split}
 & \left| S_{k-1}\left(q_1,r_1,j,\alpha_1\right)
 \right|^2\\
 \ll & \Bigg|\sum\limits_{\substack{\alpha_2\in \mathbb{Z}[i]\\ \mathcal{N}(\alpha_2)\le Q_0^{1+\varepsilon}}} \sum\limits_{q \in \mathbb{Z}[i]} 
 \Psi_2\left(\frac{q^k}{Q_0^{k/2}}\right) \cdot \Psi_2\left(\frac{(\alpha_1 + q)^k}{Q_0^{{k/2}}}\right)
 \Psi_2\left(\frac{(\alpha_2 + q)^k}{Q_0^{k/2}}\right)\times\\ & \Psi_2\left(\frac{(\alpha_1 + \alpha_2 + q)^k}{Q_0^{{k/2}}}\right)\cdot
  e\left(\Re\left(\frac{jr_1}{q_1^k}\cdot P_{k-2,\alpha_1,\alpha_2}(q)\right)\right)\Bigg|,
 \end{split}
 \end{equation*}
where
\begin{equation*}
P_{k-2,\alpha_1,\alpha_2}(q) = k(k-1)\alpha_1 \alpha_2 q^{k-2} + \cdots  
\end{equation*} 
is a polynomial of degree $k-2$ in $q$. We continue this process of repeated use of Cauchy-Schwarz and differencing until we have reached
a polynomial of degree 1. Eventually, after combining all inequalities obtained in this way, we get
\begin{equation} \label{again}
\begin{split}
& \left|S_k\left(q_1,r_1,j\right) \right|^{\kappa} \ll  Q_0^{\kappa-k+ \varepsilon}\times\\ & \sum\limits_{\substack{\alpha \in \mathbb{Z}[i]^k\\ 
\mathcal{N}\left(\alpha_1\right),...,\mathcal{N}\left(\alpha_{k-1}\right)\le Q_0^{1+\varepsilon}}} \Bigg|
\sum\limits_{q \in \mathbb{Z}[i]} \prod_{u\in \{0,1\}^{k-1}}
\Psi_2\Bigg(\frac{(u\cdot \alpha+q)^k}{Q_0^{k/2}}\Bigg) \cdot  
e\left(\Re\left(\frac{jr_1}{q_1^k}\cdot P_{1,\alpha}(q)\right)\right)\Bigg|,
 \end{split}
\end{equation}  
where we write $\kappa=2^{k-1}$, $\alpha=\left(\alpha_1,...,\alpha_{k-1}\right)$ and $u=\left(u_1,...,u_k\right)$, $u\cdot \alpha$ is the standard inner product, and $P_{1,\alpha}(q)$ takes the form
$$
P_{1,\alpha}(q) = k! \alpha_1\cdots \alpha_{k-1}\cdot \left(q + \frac{1}{2}\cdot \left(\alpha_1+\cdots +\alpha_{k-1}\right)\right).
$$

\subsection{Poisson summation}
Again, we shall apply Proposition \ref{poissongen} with $d=2$ to transform the inner-most over $q$ on the right-hand side of \eqref{again}. For 
$z=(z_1,z_2)\in \mathbb{R}^2$, we set 
$$
g(z):=\prod_{u\in \{0,1\}^{k-1}}
\Psi_2\left(\left(z_1+iz_2 + \frac{u\cdot \alpha}{\sqrt{Q_0}}\right)^k\right).
$$
Then using \eqref{obs} with 
$$
a=q \quad \mbox{and} \quad b=\overline{\frac{k! \alpha_1\cdots \alpha_{k-1}jr_1}{q_1^k}},
$$
we deduce that
\begin{equation} \label{backtor3}
\begin{split}
& \Bigg| \sum\limits_{q\in \mathbb{Z}[i]} \prod_{u\in \{0,1\}^{k-1}}
\Psi_2\Bigg(\frac{(u\cdot \alpha+q)^k}{Q_0^{k/2}}\Bigg) \cdot  
e\left(\Re\left(\frac{jr_1}{q_1^k}\cdot P_{1,\alpha}(q)\right)\right) \Bigg| \\
= & \Bigg|\sum\limits_{x\in \mathbb{Z}^2} e\left(\overrightarrow{b}\cdot x\right)g\left(\frac{x}{\sqrt{Q_0}}\right)\Bigg|.
\end{split}
\end{equation}
Now applying Proposition \ref{poissongen} with 
$B:=1/\sqrt{Q_0}$ and $f:=\tilde{g}$, the inverse Fourier transform of $g$, to the right-hand side of \eqref{backtor2}, we get
\begin{equation*}
\sum\limits_{x\in \mathbb{Z}^2} e\left(\overrightarrow{b}\cdot x\right)g\left(\frac{x}{\sqrt{Q_0}}\right) 
= Q_0\cdot \sum\limits_{y\in \overrightarrow{b}+\mathbb{Z}^2} \tilde{g}\left(\sqrt{Q_0}y\right) =
Q_0\cdot \sum\limits_{y\in -\overrightarrow{b}+\mathbb{Z}^2} \hat{g}\left(\sqrt{Q_0}y\right).
\end{equation*}
It follows that
\begin{equation} \label{newpoiss1}
\begin{split}
& \left|S_k\left(q_1,r_1,j\right) \right|^{\kappa} \ll  Q_0^{\kappa-k+ \varepsilon}\times\\ & 
\sum\limits_{\substack{\alpha \in \mathbb{Z}[i]^k\\ 
\mathcal{N}\left(\alpha_1\right),...,\mathcal{N}\left(\alpha_{k-1}\right)\le Q_0^{1+\varepsilon}}} 
\sum\limits_{\beta\in \mathbb{Z}[i]} \hat{g}\left(\sqrt{Q_0}\cdot 
\overrightarrow{\left(\beta-\overline{\frac{k!\alpha_1\cdots \alpha_{k-1}j r_1}{q_1^k}}\right)}\right).
\end{split}
\end{equation}

Here we set 
$$
\Phi_2(t):=\exp\left(-\frac{\pi}{\kappa}\cdot \sqrt[k]{|t|}\right)
$$
so that 
$$
\Psi_2(z)=\Phi_2(\mathcal{N}(z))=\exp\left(-\frac{\pi}{\kappa}\cdot \sqrt[k]{\mathcal{N}(z)}\right).
$$
It follows that
$$
g(z)=\exp\left(-\frac{\pi}{\kappa}\cdot \sum\limits_{u\in \{0,1\}^{k-1}} \left(
\left(z_1+\frac{u\cdot \alpha^{(1)}}{\sqrt{Q_0}}\right)^2+ \left(z_2 + \frac{u\cdot \alpha^{(2)}}{\sqrt{Q_0}}
\right)^2\right)\right),
$$
where 
$$
\alpha^{(1)}:=\left(\alpha_1^{(1)},...,\alpha_{k-1}^{(1)}\right):=\left(\Re(\alpha_1),...,\Re(\alpha_{k-1})\right)
$$ and 
$$
\alpha^{(2)}:=\left(\alpha_1^{(2)},...,\alpha_{k-1}^{(2)}\right):=
\left(\Im(\alpha_1),...,\Im(\alpha_{k-1})\right). 
$$
Completing the squares, it follows that
$$
g(z)=\exp\left(-\frac{\pi}{4Q_0}\cdot \sum\limits_{i=1}^2 \sum\limits_{v=1}^{k-1}
\left(\alpha_k^{(i)}\right)^2\right)\cdot
\exp\left(-\pi\cdot \sum\limits_{i=1}^2\left(z_i+\frac{\sum\limits_{v=1}^{k-1}\alpha_v^{(i)}}{2\sqrt{Q_0}}\right)^2\right).
$$
The Fourier transform of this function satisfies
\begin{equation} \label{gfourier1}
\begin{split}
\hat{g}(z)= & \exp\left(-\frac{\pi}{4Q_0}\cdot \sum\limits_{i=1}^2 \sum\limits_{v=1}^{k-1}
\left(\alpha_k^{(i)}\right)^2\right)\cdot e\left(-\frac{\sum\limits_{i=1}^2 z_i\sum\limits_{v=1}^{k-1}\alpha_v^{(i)}}{2\sqrt{Q_0}}\right)
\cdot 
\exp\left(-\pi \left(z_1^2+z_2^2\right)\right)\\
\ll & \exp\left(-\pi \mathcal{N}(z_1+iz_2)\right).
\end{split}
\end{equation}
Plugging \eqref{gfourier1} into \eqref{newpoiss1}, we get
\begin{equation} \label{finalafterps1}
\begin{split}
& \left|S_k\left(q_1,r_1,j\right) \right|^{\kappa} \ll  Q_0^{\kappa-k+1+\varepsilon}\times\\ &  
\sum\limits_{\substack{\alpha \in \mathbb{Z}[i]^k\\
\mathcal{N}\left(\alpha_1\right),...,\mathcal{N}\left(\alpha_{k-1}\right)\le Q_0^{1+\varepsilon}}} 
\sum\limits_{\beta\in \mathbb{Z}[i]} \exp\left(-\pi Q_0\mathcal{N}
\left(\beta-\overline{\frac{k!\alpha_1\cdots \alpha_{k-1}j r_1}{q_1^k}}\right)\right).
\end{split}
\end{equation}
 
\subsection{Counting}
Now we want to bound the term in the maximum in \eqref{generalKbound}. We choose $\hat{\Psi}_1$ as in \eqref{psi1}. Then
\begin{equation} \label{countbeg}
\begin{split}
& \sum\limits_{j\in \mathbb{Z}[i]\textbackslash \{0\}} \hat\Psi_1\left(\frac{\sqrt{2}jQ_0^{k/2}}{\sqrt{N}}\right) \cdot 
\left| S_k\left(q_1,r_1,j\right)\right| \ll 1+ \sum\limits_{\substack{j\in \mathbb{Z}[i]\textbackslash \{0\}\\ 
\mathcal{N}(j)\le NQ_0^{\varepsilon-k}}}  
\left| S_k\left(q_1,r_1,j\right)\right|\\
\ll & 1+ \left(\frac{N}{Q_0^{k-\varepsilon}}\right)^{1-1/\kappa} \left(\sum\limits_{\substack{j\in \mathbb{Z}[i]\textbackslash \{0\}\\ 
\mathcal{N}(j)\le NQ_0^{\varepsilon-k}}}  
\left| S_k\left(q_1,r_1,j\right)\right|^{\kappa}\right)^{1/\kappa},
\end{split}
\end{equation}
where the second line follows from  H\"older's inequality. Using \eqref{finalafterps1} and taking into account that
the contributions of $\beta$'s with 
$$
\mathcal{N}\left(\beta-\overline{\frac{k!\alpha_1\cdots \alpha_{k-1}j r_1}{q_1^k}}\right) > Q_0^{\varepsilon-1}
$$
is negligible, we deduce that
\begin{equation} \label{super1}
\begin{split}
& \sum\limits_{\substack{j\in \mathbb{Z}[i]\textbackslash \{0\}\\ 
\mathcal{N}(j)\le NQ_0^{\varepsilon-k}}}  
\left| S_k\left(q_1,r_1,j\right)\right|^{\kappa}\\ \ll & 
Q_0^{\kappa-k+1+ \varepsilon}\cdot  
\sum\limits_{\substack{j\in \mathbb{Z}[i]\textbackslash \{0\}\\ 
\mathcal{N}(j)\le NQ_0^{\varepsilon-k}}}  
\sum\limits_{\substack{\alpha \in \mathbb{Z}[i]^k\\
\mathcal{N}\left(\alpha_1\right),...,\mathcal{N}\left(\alpha_{k-1}\right)\le Q_0^{1+\varepsilon}}} 
\sum\limits_{\substack{\beta\in \mathbb{Z}[i] \\ 
\mathcal{N}\left(\beta-\overline{k!\alpha_1\cdots \alpha_{k-1}j r_1/q_1^k}\right)\le Q_0^{\varepsilon-1}}} 1\\
\ll & Q_0^{\kappa-k+1 \varepsilon}\cdot\sum\limits_{\substack{j\in \mathbb{Z}[i]\textbackslash \{0\}\\\mathcal{N}(j)\le NQ_0^{\varepsilon-k}}}  
\sum\limits_{\substack{\alpha \in \mathbb{Z}[i]^k\\ 
\mathcal{N}\left(\alpha_1\right),...,\mathcal{N}\left(\alpha_{k-1}\right)\le Q_0^{1+\varepsilon}\\ 
\left|\left|k!\alpha_1\cdots 
\alpha_{k-1}j r_1//q_1^k\right|\right|\le Q_0^{\varepsilon-1/2}}} 1\\
\ll & (NQ_0)^{(k+1)\varepsilon} \cdot Q_0^{\kappa-k+1}\cdot \Bigg(\frac{N}{Q_0^{2}}+ \sum\limits_{\substack{d\in \mathbb{Z}[i]\setminus\{0\}\\ 
\mathcal{N}(d)\le k!^2NQ_0^{k\varepsilon-1}\\   
\left|\left|dr_1/q_1^k\right|\right|\le Q_0^{\varepsilon-1/2}}} 1\Bigg)\\
\ll & (NQ_0)^{(k+1)\varepsilon} \cdot \Big(NQ_0^{\kappa-k-1}+ Q_0^{\kappa-k+1}\cdot 
\sum\limits_{\substack{l\in \mathbb{Z}[i]\\ \left|l/q_1^k\right|\le Q_0^{\varepsilon-1/2}}} 
\sum\limits_{\substack{\mathcal{N}(d)\le k!^2NQ_0^{k\varepsilon-1}\\ d\equiv l\overline{r_1} \bmod{q_1^k}}} 1\Bigg),
\end{split}
\end{equation}
where we recall that $||z||$ is the distance of $z\in \mathbb{C}$ to the nearest Gaussian integer. In the above, we have set 
$d=k!\alpha_1\cdots \alpha_{k-1}j$ if $\alpha_1,...,\alpha_{k-1}\not=0$. 

The number of residue classes modulo $q_1^{k}$ is $\mathcal{N}(q_1^k)\le Q_0^k$, and hence
\begin{equation} \label{res1}
\sum\limits_{\substack{\mathcal{N}(d)\le k!^2NQ_0^{k\varepsilon-1}\\ d\equiv l\overline{r_1} \bmod{q_1^k}}} 1 \ll 
1+\frac{N}{Q_0^{k+1-k\varepsilon}}.
\end{equation}
Further
\begin{equation} \label{lat1}
\sum\limits_{\substack{l\in \mathbb{Z}[i]\\ |l/q_1^k|\le Q_0^{\varepsilon-1/2}}} 1 \le 
\sum\limits_{\substack{l\in \mathbb{Z}[i]\\ |l|\le Q_0^{\varepsilon+(k-1)/2}}} 1\ll Q_0^{k-1+2\varepsilon}.
\end{equation}
Combining \eqref{super1}, \eqref{res1} and \eqref{lat1}, we obtain
\begin{equation} \label{comb3}
\sum\limits_{\substack{j\in \mathbb{Z}[i]\textbackslash \{0\}\\ 
\mathcal{N}(j)\le NQ_0^{\varepsilon-k}}}  
\left| S_k\left(q_1,r_1,j\right)\right|^{\kappa}
\ll (Q_0N)^{(2k+3)\varepsilon}\left(Q_0^{\kappa}+ NQ_0^{\kappa-k-1}\right),
\end{equation}
and combining \eqref{generalKbound}, \eqref{countbeg} and \eqref{comb3}, and changing $\varepsilon$ suitably, we arrive at
\begin{equation} \label{comb4}
K\ll \frac{Q_0^k}{N}+(Q_0N)^{\varepsilon}\left(\frac{Q_0^{1+k/\kappa}}{N^{1/\kappa}}+Q_0^{1-1/\kappa}\right).
\end{equation}

Now the statement in Theorem \ref{powermodZi} follows immediately from Corollary \ref{preform} and 
\eqref{comb4} after dividing the moduli into dyadic intervals and replacing $Q_0$ by $Q$.


\begin{thebibliography}{cccc}
\bibitem{Bai} S. Baier, {\it On the large sieve with sparse sets of moduli}, J. Ramanujan Math. Soc. 21 (2006), no. 3, 279--295.
\bibitem{Ba2} S. Baier, {\it The large sieve with square norm moduli in $Z[i]$}, to appear in  J. Theor. Nombr. Bordx., arXiv:1511.02470.
\bibitem{BZ1} S. Baier; L. Zhao, {\it Large sieve inequality with characters for powerful moduli}, 
Int. J. Number Theory 1 (2005), no. 2, 265--279. 
\bibitem{BZ2} S. Baier; L. Zhao, {\it An improvement for the large sieve for square moduli}, J. Number Theory 128 (2008), no. 1, 154--174.
\bibitem{BPS} W.D. Banks; F. Pappalardi; I.E. Shparlinski, {\it On group structures realized by elliptic curves over arbitrary finite fields}. 
Exp. Math. 21 (2012), no. 1, 11--25.
\bibitem{BFKS} J. Bourgain; K. Ford; S.V. Konyagin; I.E. Shparlinski, {\it On the divisibility of Fermat quotients}. 
Michigan Math. J. 59 (2010), no. 2, 313--328.
\bibitem{Hal} K. Halupczok, {\it A new bound for the large sieve inequality with power moduli}, Int. J. Number Theory 8 (2012), no. 3, 689--695.
\bibitem{HLP} G.H. Hardy; J.E. Littlewood; G. P\'olya, {\it Inequalities}. Reprint of the 1952 edition. Cambridge Mathematical Library. 
Cambridge University Press, Cambridge, 1988. 324 pp.
\bibitem{Huxl} M. Huxley, {\it The large sieve inequality for algebraic number fields}, Mathematika 15 (1968) 178-187.
\bibitem{StW} E.M. Stein; G. Weiss, {\it Introduction to Fourier analyis on Euclidean spaces}, Mir publishing house, Moscow, 1974. 336pp.
\bibitem{Zh1} L. Zhao, {\it Large sieve inequality with characters to square moduli}, Acta Arith. 112 (2004), no. 3, 297--308.
\end{thebibliography}
 \end{document}